\documentclass[11pt,reqno]{amsart}
\usepackage{amsmath,amssymb}

 \makeatletter
 \evensidemargin  \oddsidemargin
 \makeatother

\begin{document}
\thispagestyle{empty}
 \setcounter{page}{1}

\begin{center}
{\large\bf Further extension of Nadler's fixed point theorem }

\vskip.20in

{\bf  }  M. Eshaghi Gordji, M. Ramezani, H. Khodaei and H. Baghani\\
\address{Department of Mathematics,
Semnan University, P. O. Box 35195-363, Semnan, Iran}

\email{M. Eshaghi Gordji:madjid.eshaghi@gmail.com} \email{M.
Ramezani:ramezanimaryam873@gmail.com} \email{H.
Khodaei:khodaei.hamid.math@gmail.com}\email{H. Baghani:
h.baghani@gmail.com}

\vskip 5mm

\end{center}

\vskip 5mm
 \noindent{\footnotesize{\bf Abstract.}
 In this paper, we prove  a generalization of Geraghty's fixed point theorem for
 multi--valued mappings.
 \vskip.10in

  \newtheorem{df}{Definition}[section]
  \newtheorem{rk}[df]{Remark}
   \newtheorem{lem}[df]{Lemma}
   \newtheorem{thm}[df]{Theorem}
   \newtheorem{pro}[df]{Proposition}
   \newtheorem{cor}[df]{Corollary}
   \newtheorem{ex}[df]{Example}

 \setcounter{section}{0}
 \numberwithin{equation}{section}

\vskip .2in

\begin{center}
\section{Introduction}
\end{center}

Many fixed point theorems have been proved by various authors as
generalizations to Banach's contraction mapping principle. One such
generalization is due to Geraghty \cite{Ge} as follows.
\begin{thm}
Let $(X,d)$ be a complete metric space, let $f:X\to X$ be a
mapping such that for each $x,y\in X$,
$$d(f(x),f(y))\leq \alpha (d(x,y))~ d(x,y)$$
where $\alpha$ is a function from $[0,\infty)$ into $[0,1)$ which
satisfy
 the simple condition $\alpha(t_n)\to 1\Longrightarrow t_n\to
0$. Then $f$ has a fixed point $z\in X$, and $\{f^n(x)\}$ converges
to $z$, for each $x\in X$.
\end{thm}

Let $(X,d)$ be a metric space. Let $CB(X)$ denotes the collection
of all nonempty closed bounded subsets of $X$. For $A,B \in CB(X)$
and $x\in X$, define $D(x,A):=\inf\{d(x,a);a\in A\}$ and
$$ H_d(A,B):=\max\{\sup_{a\in A} D(a,B),\sup_{b\in B} D(b,A)\}.$$

 It is easy to see that $ H_d$ is a metric on $CB(X)$. $H_d$ is called the Hausdorff metric
induced by $d$.  A point $p\in X$ is said to be a fixed point of
multi--valued mapping  $T:X\to CB(X)$ if $p\in T(p).$

The fixed point theory of multi--valued contractions  was
 initiated by Nadler \cite{Na} in the following way.
 \begin{thm} $($ Nadler \cite{Na}.$)$ Let $(X,d)$ be a complete metric space and let $T$ be
 a mapping from $X$ into $CB(X)$ such that for all $x,y\in X$,
 $$ H_d(Tx,Ty)\leq r~ d(x,y)\eqno\hspace{0.5cm} (1)$$
 where, $0\leq r<1$. Then $T$ has a fixed point.
 \end{thm}
This theory was developed in different directions by many authors,
in particular, by Mizoguchi and Takahashi \cite{Mi}.
\begin{thm}(Mizoguchi and Takahashi \cite{Mi}.) Let $(X,d)$ be a complete metric space and let
$T$ be a mapping from $X$ to $CB(X)$. Assume
$$H_d(Tx,Ty)\leq \alpha (d(x,y)) ~d(x,y)\eqno\hspace{0.5cm} (2)$$
for all $x,y\in X$, where $\alpha$ is a function from $[0,\infty)$
into $[0,1)$ satisfying $\limsup_{s\to t^+} \alpha(s)<1$ for all
$t\in[0,\infty)$. Then $T$ has a fixed point.
\end{thm}
Recently, Eldred et al. \cite{El} claimed that Nadler's fixed point
theorem is equivalent to Mizoguchi and Takahashi's  fixed point
theorem. Very recently, Suzuki \cite{Su} produced an example to
disproved their claim and showed that Mizoguchi and Takahashi's
fixed point theorem is a real generalization of Nadler's theorem.

 In this paper, we extended the Geraghty's fixed point theorem to
 multi--valued mappings. Also we give an example  to show that our
 theorem is a real generalization of Nadler's.
\begin{center}
\section{Main  Result}
\end{center}

Let $S$ denotes the class of those functions $\alpha :[0,\infty)\to
[0,1)$ which satisfy the simple condition $\alpha(t_n)\to
1\Longrightarrow t_n\to 0.$\\
\begin{thm} Let $(X,d)$ be a complete metric space, let $T:X\to
CB(X)$, and suppose there exists $\alpha \in S$ such that for each
$x,y\in X$
$$H_d(Tx,Ty)\leq \alpha (d(x,y)) ~d(x,y). \eqno\hspace{0.5cm} (3)$$
Then $T$ has a fixed point.
\end{thm}
\paragraph{\bf Proof:}   Define a function $\beta$ from $[0,\infty)$ into $[0,1)$ by $$\beta
(t):=\frac{1+\alpha (t)}{2}.$$ Then the following hold:

1)$\alpha (t)<\beta (t)$ for all $t$,

2) $\beta \in S$.\\
Let $x_0\in X$ be arbitrary and fixed and let $x_1\in Tx_0$. If
$x_1=x_0$, then $x_0$ is a fixed point of $T$, and  the proof is
complete. Now, let  $x_1\neq x_0.$ Then we have
$$D(x_1,Tx_1)\leq H_d(Tx_0,Tx_1)\leq \alpha
(d(x_0,x_1))~d(x_0,x_1)<\beta (d(x_0,x_1))~d(x_0,x_1).$$ Thus
there exists $x_2\in Tx_1$ such that
$$d(x_1,x_2)\leq \beta(d(x_0,x_1))~ d(x_0,x_1).$$
Now,  if $x_1=x_2$, then $x_1$ is a fixed point of $T$, and  the
proof is complete. We suppose that $x_1\neq x_2$. Then
$$D(x_2,Tx_2)\leq H_d(Tx_1,Tx_2)\leq \alpha
(d(x_1,x_2))~d(x_1,x_2)<\beta (d(x_1,x_2))~d(x_1,x_2).$$
 Hence, there exists $x_3\in Tx_2$ satisfying
$$d(x_2,x_3)\leq \beta(d(x_1,x_2))~ d(x_1,x_2).$$
Inductively, for each positive integer number $n$, there exists
$x_{n+1}\in Tx_n$, $x_{n+1}\not= x_n$, satisfying
$$d(x_n,x_{n+1})\leq \beta(d(x_{n-1},x_n))~ d(x_{n-1},x_n).$$
To show that $\{x_n\}$ is a Cauchy sequence, we break the argument into two steps.\\\\
 Step1.  $\lim_{n\to\infty}d(x_n,x_{n+1})=0.$\\
$Proof$. Since $\beta (t)<1$ for all $t$, $\{d(x_n,x_{n+1})\}$ is
decreasing and bounded below, so
$$\lim_{n\to\infty}d(x_n,x_{n+1})=r\geq 0.$$ Assume $r>0$. Then we
have
$$\frac{d(x_{n+1},x_{n+2})}{d(x_{n},x_{n+1})}\leq \beta
(d(x_n,x_{n+1}))~~,\hspace{1cm}n=1,2,\cdots.$$ Letting $n\to
\infty$ we see that $1\leq
\lim_{n\to\infty}\beta(d(x_n,x_{n+1}))$, and since $\beta \in S$
this in turn implies $r=0$. This contradiction established Step
1.\\\\
Step 2.  $\{x_n\}$ is a Cauchy sequence.\\
$Proof.$ Assume $\displaystyle\limsup_{n,m\to\infty} d(x_n,x_m)>0$.
By triangle inequality for  positive real numbers $n,m$ and for
$y\in Tx_n$, we obtain
$$d(x_n,x_m)\leq d(x_n,x_{n+1})+d(x_{n+1},y)+d(y,x_{m+1})+d(x_{m+1},x_m).$$
This means that for every  positive real numbers $m,n$,
\begin{align*}
d(x_n,x_m)&\leq
D(x_{n+1},Tx_n)+D(x_{m+1},Tx_n)+d(x_n,x_{n+1})+d(x_m,x_{m+1})\\
&\leq H_d(Tx_m,Tx_n)+d(x_n,x_{n+1})+d(x_m,x_{m+1})\\
&\leq\beta(d(x_m,x_n)) ~d(x_n,x_m)+d(x_n,x_{n+1})+d(x_m,x_{m+1}).
\end{align*}
Hence,
$$~d(x_n,x_m)\leq (1-\beta(d(x_n,x_m)))^{-1}~(d(x_n,x_{n+1})+d(x_m,x_{m+1})).$$
Under the assumption
$~\displaystyle\limsup_{n,m\to\infty}d(x_n,x_m)>0$, Step 1 now
implies that
$$\displaystyle\limsup_{n,m\to\infty}\frac1{1-\beta
(d(x_n,x_m))}=+\infty$$
for which
$$\displaystyle\limsup_{n,m\to\infty}\beta(d(x_n,x_m))=1.$$
On the other hand we have  $\beta \in S$. It follows that
$~\displaystyle\limsup_{n,m\to\infty}d(x_n,x_m)=0~$ which is a
contradiction.\\
Now, we will complete the proof  by observing that $\{x_n\}$ is
Cauchy sequence. By completeness of  $X$, there exists $z\in X$
such that $\lim_{n\to\infty} x_n=z$. Since $T$ is continuous, then
$\lim_{n\to\infty} Tx_n=Tz$. Hence,
\begin{align*}
D(z,Tz)=D(\lim_{n\to\infty}x_{n+1},Tz)&=\lim_{n\to\infty}
D(x_{n+1},Tz)\\
&\leq\lim_{n\to\infty}H_d(Tx_n,Tz)\\
&\leq\lim_{n\to\infty}\beta(d(x_n,z))~d(x_n,z)\\
&\leq\lim_{n\to\infty}d(x_n,z)=0.
\end{align*}
On the other hand  $~Tz$ is closed. Then $z\in Tz$.
 \hfill$\Box$

 The following example shows that Theorem 2.1 is a real generalization of Nadler's.
\begin{ex}
Let $l^{\infty}$ be the Banach space consisting of all bounded real
sequences with supremum norm and let $\{e_n\}$ be a canonical basis
of $l^{\infty}$. Let $\{ \tau_n\}$ be a bounded, strictly decreasing
sequence in $(0,1)$ such that $\tau_1=\frac 12$ and for each
positive integer number $n$, $\tau_{n+1}=(1-\tau_n) \tau_n$. It is
easy to see that  $\tau_n\downarrow 0$. Put $x_n=\tau_n e_n$ and
$X_n=\{x_n, x_{n+1},\cdots\}$ for all $n\in\Bbb N$. Define a
bounded, complete subset $X$ of $l^{\infty}$ by $X=X_1$. Now define
a map $T :X\to CB(X)$ as
$$Tx_n=X_{n+1}~~,\hspace{1cm} (n\in\Bbb N)$$
and $\alpha:[0,\infty)\to[0,1)$ as

$$\alpha(t)=\begin{cases}
1-\tau_n~\hspace{1cm} & t=\tau_n ~~(n\in\Bbb
N), \\
&\qquad\qquad\\
0 & otherwise.
\end{cases}$$

Now we can see that:

(i) $T$ satisfies $(3)$ for all $x,y\in X$.

(ii) There is no $T$-invariant subset $M$ such that
$M\not=\varnothing$ and $(1)$ holds for all $x,y\in M$.

(iii) If $\alpha(t_n)\to 1$, then $t_n\to 0$.\\\\
 $proof.$ It is easy to see that the following hold:

 $\centerdot$ For $m>n$, $H(Tx_m,Tx_n)=\tau_{n+1}$.

 $\centerdot$ For $m>n$, $d(x_m,x_n)=\tau_n$.\\
 Fix $m,n\in\Bbb N$ with $m>n$, we have
 $$H(Tx_m,Tx_n)=\tau_{n+1}=(1-\tau_n) ~\tau_n=\alpha(d(x_m,x_n))~
 d(x_m,x_n).$$
It follows  $(i)$. We note that $X_n$'s are only $T$-invariant
subsets of $X$ because fix $n\in \Bbb N$, for each $k\geq n$,
$Tx_k=X_{k+1}\subseteq X_n$. If for some $k\in \Bbb N$, $(1)$ holds
for all $x,y\in X_k$, then there exists $r\in(0,1)$ such that for
each $m>n\geq k$
$$H(Tx_m,Tx_n)\leq r ~ d(x_m,x_n).$$
Hence, for all $n\geq k$ we obtain
$$\tau_n ~(1-\tau_n)\leq r~ \tau_n.$$
this shows that for each $n \geq k$, $1-r\leq \tau_n$. This is
contradiction with $\lim _{n\to\infty}\tau_n=0$. Thus, we obtain
$(ii)$. It is easy to see that $(iii)$ holds.
 On the other hand we
have
$$\limsup_{t\to 0^+}\alpha(t)=\limsup_{n\to\infty}(1-\tau_n)=1.$$
This shows that $\alpha$ does not satisfy in conditions of Mizoguchi
and Takahashi's theorem.
\end{ex}

Let $(X,d)$ be a metric space.  Let $k\in \Bbb N,$ and  $\mathcal
M^0:=X, \mathcal H^0:=d,$ for each  $i\in \{1,2,...,k\},$ put
$\mathcal M^i:=CB(\mathcal M^{i-1})$ and $\mathcal H^i:=
H_{\mathcal H^{i-1}}$. One can show that $(\mathcal M^i,\mathcal
H^i)$ is a complete metric space for all $i\in \{1,2,...,k\}$,
whenever $(X,d)$ is a complete metric space (see for example Lemma
$8.1.4$, of \cite{Rus}).
 Every mapping $T$ from $X$ into $\mathcal M^k$
is called generalized multi--valued mapping. Recently, M. Eshaghi
Gordji et al. \cite{EBKR} proved a generalized multi-valued
extension of Nadler's fixed point theorem. The question arises
here is  whether  Theorem 2.1 can be extended to generalized
multi--valued mappings
 or not?

{\small


\end{document}